\documentclass[12pt,reqno]{amsart}
\usepackage{amssymb,amsmath,amsthm}
\theoremstyle{plain} \theoremstyle{plain}
\newtheorem{theorem}{Theorem}[section]
\newtheorem{lemma}[theorem]{Lemma}

\theoremstyle{definition}

\theoremstyle{plain}
\begin{document}

\title[Lacunary Fourier series]{Lacunary Fourier series and a qualitative uncertainty
principle for compact Lie groups}

\author[E. K. Narayanan]{E. K. Narayanan}

\author[A. Sitaram]{A. Sitaram}
\address{Department of Mathematics, Indian Institute of Science,
Bangalore-560012}

\email[E. K. Narayanan]{naru@math.iisc.ernet.in} \email[A.
Sitaram]{sitaram.alladi@gmail.com}

\thanks{The first named author was supported in part by a grant
from UGC via DSA-SAP and the second named author was supported
by the Raja Ramanna Fellowship of DAE}

\begin{abstract}
We define lacunarity Fourier series on a compact connected
semisimple Lie group $G.$ If $f \in L^1(G)$ has lacunary Fourier
series, and $f$ vanishes on a non empty open subset of $G,$ then
we prove that $f$ vanishes identically. This result can be viewed
as a qualitative uncertainty principle.

\vspace*{0.1in}

\begin{flushleft}
\end{flushleft}

\vspace*{0.1in}


\end{abstract}

\maketitle

\section{Introduction}

A set $A \subset \mathbb Z$ is called $Q$-thin if it is contained
in $\mathbb N$ or $-\mathbb N$ and if $$   |n|/|m| \geq Q~~\forall
n, m \in A, |n| > |m|. $$ The set $A$ is called lacunary (in the
sense of Hadamard) if, either $A$ is empty or there are $N \in
\mathbb N$ and $Q
> 1$ such that $A \cap [N, \infty)$ and $A \cap (-\infty, -N]$ are
$Q$-thin. Then the following result is well known .

\begin{theorem}\label{the0}
Let  $A$ be a finite union of lacunary sets and $g \in L^2(\mathbb
T)$ be such that $\{m :~ \hat{g}(m) \neq 0 \} \subset A,$ then $g$
vanishes identically provided $g$ vanishes on a set of positive
measure.

\end{theorem}

This can be proved using the fact that finite union of lacunary
sets is sparse. See \cite{HJ}, page 109-110. See also \cite{Z},
Chapter V, Theorem 6.13.

Our aim in this paper is to define lacunarity on a compact
connected semisimple Lie group $G,$ using the canonical
parametrization of irreducible unitary representations of $G.$
Using our definition, we prove that, if $f \in L^1(G)$ has
lacunary Fourier series and $f$ vanishes on a non empty open
subset of $G,$ then $f$ vanishes identically. Before we state and
prove our main result, let us recall some facts.

 Let $G$ be a compact connected semisimple Lie
group and $T$ a maximal torus of $G.$ The Lie algebras of $G$ and
$T$ will be denoted by $\mathfrak g$ and $\mathfrak t,$
respectively. Let $\mathfrak g_{\mathbb C}$ and $\mathfrak
t_{\mathbb C}$ denote the complexifications of $\mathfrak g$ and
$\frak t$ respectively.

\vskip.15in \noindent We now describe the
 Weyl character formula for $G$ which will be needed later.
 For any unexplained terminology we refer to \cite{V}.

\vskip.20in
\begin{theorem}\label{the2}

Let $G$ be a compact connected semisimple Lie group and let
$R^{+}$ be a positive system of roots. Then the irreducible
characters of $G$ are in one-one correspondence with the positive
characters of $T$ (whose derivatives satisfy conditions (I) and
(II) listed below).  If $\chi$ is one such character of $T,$ the
corresponding irreducible character $ \Theta_\chi$ of $G$
(restricted to $T$) is given by
$$ \Theta_\chi = \left ( \sum_w \varepsilon(w) \chi^w \right)
\Biggr / \Delta^{+} ,$$ where $\Delta^{+}  = \sum_w \varepsilon(w)
\xi_\rho^w .$
\end{theorem}

See \cite{V}, page 40. Here $\rho$ is the half sum of positive
roots and $\xi_\rho^w = w.~\xi_\rho$ where $\xi_\rho (\textup{exp}
H) = e^{\rho( H)}$ and $w. \chi$ is the action of the Weyl group
on the characters $\chi$ of $T.$ Notice that, $\Delta^{+}(t)~
\Theta_\chi (t)$ is a finite Fourier series on $T.$

 The unitary dual of $G,$ denoted by $\hat{G},$ is described by the theorem of highest weight and is
 parameterized by a certain lattice in
 $\mathfrak t_{\mathbb C}^{*}.$ If  $\lambda \in \mathfrak
t_{\mathbb C^{*}},$ then $\lambda$ is a highest weight and gives
rise to a unitary irreducible representation of $G$ if and only if
the following two conditions are satisfied:

\vskip.10in
 \noindent I)
Positivity and integrality conditions: $\frac{2 \lambda
(H_\alpha)}{ \alpha (H_\alpha)} \in \mathbb Z^{+}$ for all $\alpha
\in R^{+},$ the set of positive roots,

\vskip.15in
 \noindent II) $\lambda$ arises from a character of $T,$ {{\em i.e.}}
$\lambda$ restricted to $\mathfrak{t}$ is the derivative of a
character on $T.$

\vskip.10in

 \noindent
 It follows that $\hat{G}$ can be identified with a subset of a (full) lattice contained in the
dual of $T.$ For more on compact semsimiple Lie groups we refer to
\cite{K} and \cite{V}.

Next, we define lacunary Fourier series on $G.$ Let $f$ be an
$L^1$ function on $G$ and $\pi(f)$ denote the Fourier transform of
$f$ defined by $$ \pi(f) = \int_G~f(x)~\pi(x)~dx,$$ where $\pi$ is
an irreducible unitary representation of $G.$ As mentioned earlier
$\hat{G}$ may be identified with a subset of a lattice in the dual
of the maximal torus $T.$ This is done as follows: We identify the
lattice generated by the elements in $\mathfrak t_{\mathbb C}^{*}$
satisfying (I) and (II) above with a subset of the lattice
$\mathbb Z^k$ where
 $k = \textup{dim}~T.$ We say that, $f \in L^1(G)$ has lacunary Fourier series if the set
$E = \{ \pi:~ \pi(f) \neq 0 \}$ (considered as the lattice defined
earlier) satisfies the following:

\begin{equation}\label{eqn1}
 \bigcup_{w \in W}~ w. E ~\subset E_1 \times E_2 \times  \cdots \times
E_k
\end{equation} where each $E_j \subset \mathbb Z,~j = 1, 2, \cdots k$ is a finite union of
lacunary sets.  {{\it We emphasize that our definition of
lacunarity differs from that of Helgason}} (see \cite{H}).

\section {Proof of the main result}

In this section we prove the main result. We start with the
following lemma.

\vskip.20in

\begin{lemma}\label{lem1} Let $g \in C(\mathbb T^k)$ and $spec(g) =
\{ \hat{g}( {\bf{m}}) \neq 0:~ {\bf{m}} \in \mathbb Z^k \} \subset
E_1 \times E_2 \times \cdots  \times E_k,$ where each $E_j \subset
\mathbb Z, ~j=1, 2, \cdots k$ is a finite union of lacunary sets.
If $g$ vanishes on a non empty open subset of $\mathbb T^k,$ then
$g$ vanishes identically.
\end{lemma}

\begin{proof}
This follows from Theorem \ref{the0} and  induction on $k.$
\end{proof}

 \vskip.15in \noindent
 Now we are in a position to state and prove the
main result of this paper.

 \vskip.2in
\begin{theorem}\label{the1}
Suppose $f \in L^1(G)$ has lacunary Fourier series. If $f$
vanishes on a non empty  open subset of $G,$  then $f \equiv 0.$

\end{theorem}

\begin{proof} Let $spec(f) \subset \mathbb Z^k$ correspond to
$\{\pi \in \hat{G}: \pi(f) \neq 0\}.$ Notice that if $f$ has
lacunary Fourier series then so do the translates of $f.$ Hence,
without loss of generality we may assume that $f$ vanishes in a
neighborhood of the identity in $G.$ Convolving $f$ with smooth
approximate identities supported in small neighborhoods of the
identity of $G,$ we may assume that $f$ is continuous, indeed
smooth.

 Now, consider the function $F_f$
defined as
$$F_f (x) = \int_G~f(g x g^{-1})~dg.$$
Since one can choose a metric on $G$ which is biinvariant, it is
clear that $F_f$ vanishes in a small enough neighborhood of
identity in $G.$ Also, $F_f$ is a central function and so is
determined by its restriction (still denoted by $F_f$) to the
maximal torus $T.$ A simple computation using Schur's
orthogonality relations shows that
$$ F_f(t) = \sum _{\pi \in \hat{G}}~\textup{Trace}(\pi(f))~\Theta_\pi (t),$$
where $t \in T$ and $\Theta_\pi$ is the character of $\pi$ defined
by $\Theta_{\pi}(x) = \textup{Trace}(\pi(x)).$  By Weyl's
character formula we get, $$ \Delta^+ (t) F_f(t) = \sum_{\pi \in
\hat{G} } ~\textup{Trace} (\pi(f)) ~\left ( \sum_w \varepsilon(w)
\chi_\pi^w (t) \right),$$ where the positive character $\chi_\pi$
of $T$ corresponds to $\pi \in \hat{G}.$  Notice that, the set of
 nonzero Fourier coefficients of $\Delta^{+}(t)~ F_f (t)$ is
contained in $\bigcup_{w \in W} ~w. spec(f)$ which satisfies
(\ref{eqn1}). Since $F_f$ vanishes on an open set, by Lemma
\ref{lem1}, $F_f$ vanishes identically. It follows that
$\textup{Trace}(\pi(f))$ is zero for $\pi \in \hat{G}.$

Next, notice that any translate of $f$ has lacunary Fourier
series. If $g$ varies in a small enough neighborhood of the
identity in $G,$ then applying the above argument to the
translated function $^g f (x) = f(gx)$ we obtain that $
\textup{Trace} (\pi(g^{-1}) \pi(f))$ is zero for all $g$ in a
small neighborhood and all $\pi \in \hat{G}.$ For each $\pi,$  $g
\rightarrow \textup{Trace}(\pi(g^{-1}) \pi(f))$ is a real analytic
function and hence it follows that $\textup{Trace} (\pi(g^{-1})
\pi(f))$ is identically zero for each $\pi.$ This clearly implies
that $f$ is zero, by the inversion formula.
\end{proof}

\vskip.10in \noindent {\bf Remark I:}~ Although we have assumed
semisimplicity, we can modify the proof so that the theorem holds
for all compact connected Lie groups.

\vskip.10in \noindent{\bf Remark II:}~If $G = U(n),$ then
$\hat{G}$ can be canonically identified with a subset of $\mathbb
Z^n$ (see \cite{V}). The Weyl group $W$ is the permutation group
$S_n.$ If $A = A_1 \times A_2 \times \cdots A_n$ where each $A_j
\subset \mathbb Z, ~j =1, 2, \cdots n,$ is lacunary, then it is
easy to see that $A$ satisfies (\ref{eqn1}). (Note that $U(n)$ is
not semisimple.)

\vskip.10in \noindent {\bf Remark III:}~Ideally one would like to
get the exact analogue of Theorem \ref{the0}, but the averaging
method used here to convert the problem to one on the torus works
only if we assume that the set on which $f$ vanishes is open. If
$f$ vanishes on a set of positive measure, then $F_f$ considered
as a function on $\mathbb T$ need not necessarily vanish on a set
of positive measure on $\mathbb T.$ However, if $f$ vanishes on a
neighborhood of the identity in $G,$ then $F_f$ considered as a
function on $\mathbb T$ vanishes on a neighborhood of the identity
in $\mathbb T.$

\section{Uncertainty principles}

To quote G. B. Folland, " the uncertainty principle is a meta
theorem in harmonic analysis that can be summarized as follows:
{{\it  A non zero function and its Fourier transform cannot both
be sharply localized}}.''  For more on uncertainty principles on
compact groups we refer to \cite{PS} and \cite{GK} and for
uncertainty principles in general to \cite{FS}.

A function on $\mathbb Z^n$ can be thought of as localized, if it
is supported on a lacunary set. The opposite of a localized
function is a spread out function. A function whose closed support
is all of $G$ can be thought of as a highly spread out function.
Thus Theorem \ref{the1} says that if the Fourier transform of a
function on $G$ is localized, then the function itself has to be
necessarily spread out---thus illustrating the uncertainty
principle in the context of compact groups.

 \vskip.20in \noindent {{\bf Acknowledgement:}}~We thank
G. B. Folland for drawing our attention to Chapter V, Theorem 6.13
of \cite{Z}.

\end{document}